\documentclass[a4paper]{article}
\usepackage[a4paper,left=3cm,right=3cm]{geometry}
\usepackage[T1]{fontenc}

\usepackage{amsmath,amssymb,amsthm}

\usepackage{hyperref}
\providecommand{\email}[1]{\href{mailto:#1}{#1}}

\newtheorem{lemma}{Lemma}

\title{Multi-Objective Mixed Integer Programming: An Objective Space Algorithm}

\author{William Pettersson \\
School of Computing Science, University of Glasgow \\ Glasgow G12 8QQ, \textsc{United Kingdom}
\email{william@ewpettersson.se}
\and Melih Ozlen \\
School of Science, RMIT University, \\ Victoria~3000, \textsc{Australia}
\email{melih.ozlen@rmit.edu.au}
}

\date{15 December 2018}

\begin{document}
\maketitle

\begin{abstract}
  This paper introduces the first objective space algorithm which can exactly find all supported and non-supported non-dominated solutions to a mixed-integer multi-objective linear program with an arbitrary number of objective functions. This algorithm is presented in three phases. First it builds up a super-set which contains the Pareto front.
  This super-set is then modified to not contain any intersecting polytopes.
  Once this is achieved, the algorithm efficiently calculates which portions of the super-set are not part of the Pareto front and removes them, leaving exactly the Pareto front.
\end{abstract}

\section{Introduction}
Many problems that arise in the real world are mixed-integer problems, problems that have both discrete and continuous variables. For instance, the facility location problem can be formulated using continuous variables the shipment amounts and discrete variables to describe the selection of the sites. Vehicle routing problems can also utilise continuous variables for loads carried, and can use discrete variables to represent vehicles. For this paper, we will only consider linear objective functions, and as such will refer to mixed-integer linear programming or MILP problems.

These real world problems often also have several objectives (such as maximising profit or reducing safety risks or environmental harm) which often conflict with each other. Algorithms for such multi-objective problems will often present a set of non-dominated solution points. This can highlight the trade-offs between the various objectives to a decision maker, allowing for a more informed decision.
The resulting problems are then multi-objective mixed-integer linear programming problems, or MOMILP problems.

Research into MOMILPs is relatively newer than other multi-objective programming problems, and it started with a special case including only binary (0 or 1) variables, Multi-objective Mixed Binary Linear Programming (MOMBLPs). A branch-and-bound approach to MOMBLPs was introduced by Mavrotas and Diakoulaki~\cite{MAVROTAS1998530} and has since been improved~\cite{MAVROTAS200553, VINCENT2013498}. We note that their original algorithm finds only supported non-dominated points --- points which lie on some non-dominated facet of the objective space. Vincent et.~al.~discuss this in Section 3.1 of their paper~\cite{VINCENT2013498}, and give an explicit algorithm to find all non-dominated solutions to bi-objective mixed-integer linear programs (BOMBLPs).
Further work on the bi-objective algorithms include a local branch-and-bound algorithm~\cite{Belotti2013Slices} and an exact bi-objective algorithm based on the $\varepsilon$- and Tabu-constraint methods~\cite{SOYLU2016Exact}.

Our new algorithm utilises existing multi-objective integer programming (MOIP) techniques, combined with slices~\cite{Belotti2013Slices} and Hamming constraints~\cite{SOYLU2016Exact}, to recursively build up a set which contains the complete non-dominated set of solutions.
We then use techniques similar to those from~\cite{VINCENT2013498}, but generalised to arbitrary dimensionality, to reduce this set to the exact set of non-dominated solutions.

%Section \ref{sec:background} gives a background into both optimisation theory and techniques, as well as a brief overview of the polytope theory we utilise.
%Section \ref{sec:algo} goes over our new algorithm in some detail, and Section \ref{sec:discussion} gives some concluding discussion and remarks.

%\section{Background}\label{sec:background}

\section{Optimisation background}

A multi-objective programming problem is defined as
\begin{align*}
  \min \; & \boldsymbol{f} = (f_1(\boldsymbol{x}), f_2(\boldsymbol{x}), \ldots, f_k(\boldsymbol{x})) \\
  \text{ s.t. } & A\boldsymbol{x} \leq B \\
  & \boldsymbol{x} \in X
\end{align*}
where $X$ is a specified search space.
It is not clear that a single optimal solution to such a problem should exist, instead, often one searches for non-dominated solutions. We say a solution
%$(f_1(x_1), \ldots, f_k(x_1))$
$\boldsymbol{f}(\boldsymbol{x_1})$
dominates a solution $\boldsymbol{f}(\boldsymbol{x_2})$
%$(f_1(x_2), \ldots, f_k(x_2))$
if $f_i(\boldsymbol{x_1}) \leq f_i(\boldsymbol{x_2})$ for all $i$ and there is some $j$ such that $f_j(\boldsymbol{x_1}) < f_j(\boldsymbol{x_2})$. The non-dominated solutions are then exactlx those solutions which are never dominated. Together, these are called the Pareto front.

Exact algorithms for both the pure integer ($X = \mathbb{Z}^n$) and pure continuous ($X = \mathbb{R}^n$) variants have been well studied, however we give an exact algorithm to solve multi-objective mixed-integer programming (MOMIP) problems (where $X=\mathbb{Z}^n \times \mathbb{R}^m$).

Belotti et.~al.~\cite{Belotti2013Slices} introduce the notion of a ``slice'' of a MOMIP, which is also used in an exact bi-objective algorithm~\cite{SOYLU2016Exact}. Given a MOMIP where $X=\mathbb{Z}^n \times \mathbb{R}^m$, a ``slice'' of this problem is created by fixing the variables $(x_1,x_2,\ldots,x_n) \in\mathbb{Z}^n$. Since all of the integer variables are now fixed, only continuous variables remain and we are left with a regular linear program. This allows the use of various techniques from linear programming which would not apply in a mixed-integer setting.

Given vectors $\boldsymbol{x} = (x_1, x_2, \ldots, x_n)$ and $\boldsymbol{\hat x} = (\hat x_1, \hat x_2, \ldots, \hat x_n)$, the Hamming distance $H(\boldsymbol{x}, \boldsymbol{\hat x})$ is defined as the number of coordinates in which $\boldsymbol{x}$ and $\boldsymbol{\hat x}$ differ.
More formally, $H(\boldsymbol{x}, \boldsymbol{\hat x}) = |\{k | x_k \neq \hat x_k\}|$.
Soylu and Yildiz~\cite{SOYLU2016Exact}, introduce the notion of a ``Hamming constraint'', $H(\boldsymbol{x}, \boldsymbol{\hat x}) \geq 1$, and also give a linearisation. We will use Hamming constraints to avoid solutions in a given slice.

\section{Polytopes}
A polytope is the generalisation of a convex polygon to higher dimensions.
A simplex (plural simplices) is a polytope whose vertices are affinely independent.
Much of the complexity that arises in this algorithm arises due to the geometry of working in Euclidean space in an arbitrary number of dimensions, specifically when it comes to intersections and set differences of various polytopes.
A polytope can be defined as the convex hull of a set of vertices in some Euclidean space.
This is called the V-representation of the polytope.
A polytope can also be defined as the intersection of a set of halfspaces of Euclidean space, generally where each halfspace corresponds to a facet of the polytope.
This representation will be called the H-representation.
Converting between H-representations and V-representations is not computationally trivial, but is well studied in its own right. The more common methods include the Double Description Method~\cite{Fukuda1996DD}, and Fourier Motzkin elimination~\cite{DANTZIG1973288}.
For further reading on the details of either of these, or polytopes in general, we guide the reader to Lectures on Polytopes~\cite{Ziegler1995}.

We will be comparing polytopes in the solution space. We will say a polytope $P$ {\em completely dominates\/} (or just dominates) another polytope $P'$ if for every point $p' \in P'$, there exists some point $p\in P$ such that $p$ dominates $p'$.
A polytope is {\em completely dominated\/} if there is some polytope which dominates it, and {\em completely non-dominated\/} if there is no such polytope.
We will say a polytope $P$ {\em partially dominates\/} another polytope $P'$ if there is no point $p' \in P'$ such that $p'$ dominates some point $p\in P$, and some point $p\in P$ such that $p$ dominates some point $p'\in P'$.
Note that we do allow a $p'\in P'$ and $p\in P$ such that $p' = p$, hence the term partial.

\section{Overview of the algorithm}\label{sec:algo}

As is common in multi-objective optimisation, our algorithm operates in the solution space of the problem.
The algorithm operates in three distinct steps. To make it easier to follow, and to allow further research to improve this work, we first present an overview of each portion of the algorithm, and then explain each step on its own. 
Note that we explain here a simplified version of this algorithm. The journal version will present these algorithms in their fullness, as well as many efficiency improvements and run-time analyses.

The first phase in our algorithm builds a set of simplices $\mathcal{C}$ such that the Pareto front is contained within $\mathcal{C}$.
Note that this does not mean that each simplex in $\mathcal{C}$ is either in the Pareto front or not.
Instead, for any point $p$ within the Pareto front, there is some $P\in\mathcal{C}$ such that $p\in P$.
We obviously keep this condition throughout our new algorithm.

The second phase then breaks up some of these simplices into even smaller polytopes. This is done to ensure that after this phase of the algorithm, no two polytopes in our set have a non-empty intersection. This allows us to work with each polytope as a whole.

The third phase analyses all of these polytopes, and determines when one polytope either partially or completely dominates some other polytope. 
Rather than breaking the polytopes into yet smaller pieces, we instead mark where each polytope is not dominant.
This marking is such that in our solution each individual polytope $P$ which contains some portion of the Pareto front will be marked with its own set of polytopes $\{R_1,\ldots,R_t\}$ such that $P \setminus \{R_1,\ldots,R_t\}$ only contains portions of the Pareto front.
%We once again break our polytopes up so that we can remove the portion of the polytope which is dominated.
%Once this is done, the Pareto front to the MOMILP is exactly the non-dominated (non-marked) sections of the polytopes.

\subsection{Phase 1: Collecting The Pieces}
The first phase combines integer and linear programming techniques to find a set $\mathcal{C}$ of simplices such that the Pareto front is contained within $\mathcal{C}$.
Note that this does not mean that each simplex in $\mathcal{C}$ is either in the Pareto front or not.
Instead, for any point $y$ within the Pareto front, there is some $Y\in\mathcal{C}$ such that $y\in Y$.

The algorithm used for this phase is a modification of MOIP\_AIRA~\cite{Ozlen2014}
with Hamming constraint and slices~\cite{SOYLU2016Exact}.
Following MOIP\_AIRA, as each solution point $p$ is found, we use Benson's algorithm~\cite{LOHNE2017807} to find all polytopes which are optimal in the given slice and contain $p$.
Hamming constraints are then added to the problem to avoid the slice containing $p$ (and thus avoid detecting the same polytope again).
To ensure that any other polytopes in this slice are found, when a new point $p$ is discovered the algorithm also checks to see if $p$ is ``beyond'' a vertex on a known polytope. If so, the known vertex is checked for additional polytopes which contain it. If none are found, then the Hamming constraint corresponding to this vertex is removed.

The correctness of this approach can be shown by contradiction.
A rough outline of the proof is as follows. Assume the existence of a polytope which was not found by the algorithm, which also gives a lexicographically least vertex $y$ of this polytope.
The only constraints we add to our problem that would cause $y$ to be missed are Hamming constraints.
As $y$ was never found, any other polytope which contains $y$ was also not found. The existence of such a polytope is given by convexity of the solution of the linear program which results from restricting our original problem to this slice. This reduction leads to a polytope with a lower least vertex, and thus a contradiction.

This phase builds our set of polytopes, which we call $\mathcal{C}$.

\subsection{Phase 2: Carving Out Intersections}

For this phase, we assume we have a set $\mathcal{C}$ in objective space such that for any non-dominated solution $y$, there exists a polytope $P\in\mathcal{C}$ such that $y\in P$.
From this set we build a set of polytopes in the objective space such that each polytope is a subset of some $P\in\mathcal{C}$ and no two polytopes have a non-empty intersection.

Consider every possible intersection of polytopes in $\mathcal{C}$. If the intersection is not empty, we add it to our set. To avoid having intersecting polytopes, we must subtract the intersection from each input polytope. However this might also produce two (or more) connected components, and these components may be non-convex.

To avoid these non-convex components, we want to break up any non-convex components into convex polytopes.
We do this by intersecting each input polytope with each of the inverse of the halfspaces which define the intersection polytope.
This procedure may have to be repeated iteratively, but does terminate.

At this point, we have a set of polytopes $\mathcal{C}^*$, of which no pair has non-empty intersection. Thus it satisfies the following.

\begin{lemma}\label{lemma:intersect}
Let $\mathcal{C}^*$ be a set of regions with no two regions having a non-empty intersection. Let $P,P'\in\mathcal{C}^*$ and let $y\in P$. If $y$ is dominated by some other solution $y'\in P'$, then there is no point $z\in P$ such that $z$ dominates any point $z'\in P'$.
\end{lemma}
A complete proof will be included in the journal version of this paper. The basic idea is to consider the two lines which run from $y$ to $z$ and from $y'$ to $z'$ respectively. These lines will be contained within $P$ and $P'$ respectively as both are convex, and must intersect, giving a contradiction.

This lemma tells us that we do not need to consider any further intersections between any polytopes to find the Pareto front.
However, we may still have some polytopes which are partially dominated by others.
A partially dominated polytope is a polytope which contains a point which is dominated by a point in a second polytope.

\subsection{Phase 3: Determining Dominance}
For this phase, we have a set of polytopes in the objective space such that no two polytopes have a non-empty intersection. 
A naive approach at this point would be to further split up each polytope into smaller polytopes, such that we get a set of polytopes where each polytope is either completely non-dominated, or completely dominated.
This can, however, lead to a significant increase in the memory usage of the algorithm.

Instead, for any set of polytopes where one polytope partially dominates the others, we build exactly the polytope which describes the portion of the objective space where this dominance holds.

This is achieved by projecting each polytope in each coordinate direction, finding the intersection in each dimension and reversing the projection.
This procedure gives a new polytope $P^*$ which marks the region of objective space where all of the input polytopes have a solution. By Lemma~\ref{lemma:intersect} only one of these input polytopes can be dominant, so for any other polytope we mark that it is not dominant in $P^*$. Note that we do not add $P^*$ as a potentially dominant polytope --- it may well lie partially outside the feasible space for the problem.

The second algorithm used in Phase 3 simply iterates over all possible subsets to determine the Pareto front. The Pareto front is presented as a set of regions, with each region being a polytope with a finite (and possibly zero) number of other polytopes removed from it.
That is, these regions are not necessarily convex.

\section{Discussion}\label{sec:discussion}

We have presented here the first multi-objective mixed-integer linear programming algorithms which can exactly determine all supported and non-supported non-dominated solutions for an arbitrary number of objective functions.
There are many improvements that can be made% in terms of efficiency
in this algorithm which are obvious
and have not been
discussed.
More complete details on these improvements, as well as full run-time analyses will be given in the journal version of this paper.
We hope this paper will spur further research in the field of exact algorithms for MOMIP problems.

% Sections that will go in second font

% Acknowledgement
\section{Acknowledgements}
The first author was supported by grant EP/P028306/1 from the Engineering and Physical Sciences Research Council.

% References

\bibliographystyle{plain}
\bibliography{bib}

\begin{thebibliography}{10}

\bibitem{Belotti2013Slices}
P.~Belotti, B.~Soylu, and M.~Wiecek.
\newblock A branch-and-bound algorithm for biobjective mixed-integer programs.
\newblock 2013.

\bibitem{DANTZIG1973288}
George~B Dantzig and B~Curtis Eaves.
\newblock Fourier-motzkin elimination and its dual.
\newblock {\em Journal of Combinatorial Theory, Series A}, 14(3):288 -- 297,
  1973.

\bibitem{Fukuda1996DD}
Komei Fukuda and Alain Prodon.
\newblock Double description method revisited.
\newblock In Michel Deza, Reinhardt Euler, and Ioannis Manoussakis, editors,
  {\em Combinatorics and Computer Science}, pages 91--111, Berlin, Heidelberg,
  1996. Springer Berlin Heidelberg.

\bibitem{LOHNE2017807}
Andreas L\"{o}hne and Benjamin We{\ss}ing.
\newblock The vector linear program solver bensolve --- notes on theoretical
  background.
\newblock {\em European Journal of Operational Research}, 260(3):807 -- 813,
  2017.

\bibitem{MAVROTAS1998530}
G.~Mavrotas and D.~Diakoulaki.
\newblock A branch and bound algorithm for mixed zero-one multiple objective
  linear programming.
\newblock {\em European Journal of Operational Research}, 107(3):530 -- 541,
  1998.

\bibitem{MAVROTAS200553}
G.~Mavrotas and D.~Diakoulaki.
\newblock Multi-criteria branch and bound: A vector maximization algorithm for
  mixed 0-1 multiple objective linear programming.
\newblock {\em Applied Mathematics and Computation}, 171(1):53 -- 71, 2005.

\bibitem{Ozlen2014}
Melih Ozlen, Benjamin~A. Burton, and Cameron A.~G. MacRae.
\newblock Multi-objective integer programming: An improved recursive algorithm.
\newblock {\em Journal of Optimization Theory and Applications},
  160(2):470--482, Feb 2014.

\bibitem{SOYLU2016Exact}
Banu Soylu and Gazi~Bilal Yildiz.
\newblock An exact algorithm for biobjective mixed integer linear programming
  problems.
\newblock {\em Computers \& Operations Research}, 72:204 -- 213, 2016.

\bibitem{VINCENT2013498}
Thomas Vincent, Florian Seipp, Stefan Ruzika, Anthony Przybylski, and Xavier
  Gandibleux.
\newblock Multiple objective branch and bound for mixed 0-1 linear programming:
  Corrections and improvements for the biobjective case.
\newblock {\em Computers \& Operations Research}, 40(1):498 -- 509, 2013.

\bibitem{Ziegler1995}
G\"unter~M. Ziegler.
\newblock {\em Lectures on polytopes}, volume 152 of {\em Graduate Texts in
  Mathematics}.
\newblock Springer-Verlag, New York, 1995.

\end{thebibliography}

\end{document}